\documentclass[12pt]{article}%
\usepackage{amsfonts}
\usepackage{sw20bams}
\usepackage{amsmath}
\usepackage{amssymb}
\usepackage{graphicx}%
\providecommand{\U}[1]{\protect\rule{.1in}{.1in}}
\begin{document}

\title{Multisum Sets}
\author{Steven Finch}
\date{July 6, 2024}
\maketitle

\begin{abstract}
Complete infinite multisum sets are eventually linear. \ After 30 years of
sitting in a file cabinet, the proof (thanks to James H. Schmerl) is brought
from darkness into light.

\end{abstract}

\footnotetext{Copyright \copyright \ 2024 by Steven R. Finch. All rights
reserved.}Let $S$ denote a nonempty set of positive integers. \ If $s,t\in S$,
then $s+t$ is called a \textbf{sum}. \ The set $S$ is a \textbf{sum set} if
all of its sums are in $S$. \ Clearly such an $S$ is infinite. \ It can be
proved that $S$ is eventually linear \cite{NW-mltsum, SS-mltsum, Hg-mltsum},
i.e., there exist integers $N$ and $k$ such that for all $n>N$, $n\in S$ if
and only if $k\mid n$.

Let us now start over. \ If $s,t,u,v\in S$ satisfy $s+t=u+v$ and are distinct,
except possibly $u=v$, then $s+t$ is called a \textbf{multisum}. \ The set $S$
is a \textbf{multisum set} if all of its multisums are in $S$. \ Such sets can
be finite. \ For example, $\{1,3,7\}$ is vacuously multisum whereas
$\{1,3,5,6\}$ is non-vacuously multisum ($1+5=3+3$ is contained in the set).
\ An infinite multisum set $S$ is \textbf{complete} if every sufficiently
large element is a multisum. \ Such a set can be proved to be eventually
linear, and this task will occupy us for the remainder of the paper.

\section{Schmerl's Theorem}

Schmerl \cite{S1-mltsum} proved the following more general result. \ Since his
work has remained unpublished, as far as is known, it seems important to
record it for posterity's sake.

Let $a_{1}<a_{2}<a_{3}<\ldots$ be a sequence of positive integers and let
$n>0$ be such that:

\begin{itemize}
\item whenever $m>n$, then $a_{m}=a_{i}+a_{j}=a_{r}+a_{s}$ for some $i<r\leq
s<j$; and

\item whenever $a=a_{i}+a_{j}=a_{r}+a_{s}>a_{n}$, where $i<r<s<j$, then
$a=a_{m}$ for some $m>n$.
\end{itemize}

\noindent Then the set $A=\{a_{1},a_{2},a_{3},\ldots\}$ is eventually linear.

\section{Lemma 1}

Assume that $d,a,b,a+d,b+d$ are distinct and in $A$, and that $a_{n}<k=a+b+d$.
\ Then all multiples of $k$ are in $A$.\medskip

\textbf{Proof: }We show by induction on $m$ that $m\,k\in A$ for each $m\geq
1$. \ This is true for $m=1$ because%
\[
k=a+b+d=[a+d]+b=[b+d]+a.
\]
Suppose that $(m-1)k+d,(m-1)k+a+d,(m-1)k+b+d,m\,k\in A$. \ The following
identities suffice to prove that $(m+1)k\in A$:%
\[
m\,k+d=[(m-1)k+b+d]+[a+d]=[(m-1)k+a+d]+[b+d],
\]%
\[
m\,k+a+d=[m\,k+d]+a=m\,k+[a+d],
\]%
\[
m\,k+b+d=[m\,k+d]+b=m\,k+[b+d],
\]%
\[
(m+1)k=[m\,k+a+d]+b=[m\,k+b+d]+a.
\]

\section{Lemma 2}

Assume $d_{1},d_{2},x,y$ are such that $x\neq d_{1}\neq d_{2}\neq y$ and%
\[
\{d_{1},2d_{1},x,x+d_{1},x+2d_{1}\}\cup\{d_{2},2d_{2},y,y+d_{2},y+2d_{2}%
\}\subseteq A,
\]%
\[
\{x,x+d_{1}\}\cap\{y,y+d_{2}\}\neq\varnothing,
\]
and $d_{1}+d_{2}>a_{n}$. \ Then there is $k\geq1$ all of whose multiples are
in $A$.\medskip

\textbf{Proof:} We will use Lemma 1 to obtain the existence of $k$ by
exhibiting $a,b,d$ such that $d,a,b,a+d,b+d$ are distinct and in $A$.

Without loss of generality, assume $d_{1}<d_{2}$. \ The hypothesis
$\{x,x+d_{1}\}\cap\{y,y+d_{2}\}\neq\varnothing$ naturally leads to four cases.

\begin{enumerate}
\item[(i)] Suppose $x+d_{1}=y+d_{2}$. \ Then let $d=x+d_{1}$, $a=d_{1}$ and
$b=d_{2}$. \ For example, $d\neq b$ since $x+d_{1}\neq x+d_{1}-y=d_{2}$. \ As
another example, $b\neq a+d$ since $d_{2}=x+d_{1}-y\neq d_{1}+x+d_{1}$.

\item[(ii)] Suppose $x+d_{1}=y$. \ Then let $d=y$ and $a=d_{1}$. \ If
$d_{2}=x+2d_{1}$, then let $b=2d_{2}$; if $d_{2}\neq x+2d_{1}$, then let
$b=d_{2}$.

\item[(iii)] Suppose $x=y+d_{2}$. \ Then let $d=x$, $a=d_{1}$ and $b=d_{2}$.

\item[(iv)] Suppose $x=y$. \ Then let $d=x$ and $a=d_{1}$. \ If $d_{2}%
=d_{1}+x$, then let $b=2d_{2}$; if $d_{2}\neq d_{1}+x$, then let $b=d_{2}$.
\end{enumerate}

\section{Proof of Theorem, Part One}

Let $M=6n-4$, $I=\{a_{i}:1\leq i<M\}$ and $J=\{a_{j}:n<j\leq M\}$. \ For each
$t\in J$, let $T_{t}=\{x,y,z\}\subseteq I$ be such that for some $w\in I$, we
have $x<y\leq w<z$ and $t=x+z=y+w$. \ Let%
\[
D=\left\{  d\in I:\exists\text{ distinct }r,s,t\in J\text{ such that }d\in
T_{r}\cap T_{s}\cap T_{t}\right\}  .
\]
Now consider some $d\in D$, where $d\in T_{r}\cap T_{s}\cap T_{t}$ and
$r<s<t$. \ At least one of the following six alternatives must hold:

\begin{enumerate}
\item[(1)] $r=2d$ and $t=s+d;$

\item[(2)] $s=2d$ and $t=r+d;$

\item[(3)] $t=2d$ and $s=r+d;$

\item[(4)] $r\neq2d$, $s\neq2d$ and $s\neq r+d;$

\item[(5)] $r\neq2d$, $t\neq2d$ and $t\neq r+d;$

\item[(6)] $s\neq2d$, $t\neq2d$ and $t\neq s+d.$
\end{enumerate}

\noindent If either (4), (5) or (6)\ hold, then by making the following
choices for $a$ and $b$, respectively, the hypothesis of Lemma 1 will be
satisfied:%
\[
a=r-d\text{ and }b=s-d;
\]%
\[
a=r-d\text{ and }b=t-d;
\]%
\[
a=s-d\text{ and }b=t-d.
\]
Hence we can assume that one of (1), (2), (3)\ holds. \ If (1)\ or (3)\ holds,
then%
\[
s\neq2d\text{ and }\{2d,s-d,s,s+d\}\subseteq I\cup\{a_{M}\};
\]
and if (2)\ holds, then%
\[
r\neq2d\text{ and }\{2d,r-d,r,r+d\}\subseteq I\cup\{a_{M}\}.
\]
In any case, there is $x\neq d$ such that $\{2d,x,x+d,x+2d\}\subseteq
I\cup\{a_{M}\}$. \ Let $S_{d}=\{x,x+d\}\subseteq I$, where $x$ is as just described.

As noted, if $d\in T_{r}\cap T_{s}\cap T_{t}$, where $r<s<t$, then
$2d\in\{r,s,t\}.$ \ It follows that if $r<s<t<u$ are in $J$, then $T_{r}\cap
T_{s}\cap T_{t}\cap T_{u}=\varnothing$. \ Consequently%
\[
3\left\vert D\right\vert +2\left(  \left\vert I\right\vert -\left\vert
D\right\vert \right)  \geq3\left\vert J\right\vert ,
\]
which implies that $\left\vert D\right\vert \geq3n-2$. \ Thus%
\[
2\left\vert D\right\vert \geq6n-4>6n-5=\left\vert I\right\vert ,
\]
so there are distinct $d_{1},d_{2}\in D$ for which $S_{d_{1}}\cap S_{d_{2}%
}\neq\varnothing$. \ Let $S_{d_{1}}=\{x,x+d_{1}\}$ and $S_{d_{2}}%
=\{y,y+d_{2}\}$. \ Then $d_{1},d_{2},x,y$ are as in the hypothesis of Lemma 2,
yielding a $k$ all of whose multiples are in $A.$

\section{Proof of Theorem, Part Two}

Let $k$ be the least positive integer such that all sufficiently large
multiples of $k$ are in $A$. \ (By the preceding section, we know that such a
$k$ exists.) \ Let $M$ be such that $m\geq M$ implies $m\,k\in A$.

Suppose that not all sufficiently large elements of $A$ are multiples of $k$.
\ Then there is $r$ such that $1\leq r<k$ and there are $x,s\geq1$ for which
\ $x\equiv r$ ($\operatorname{mod}k$) and $x,x+s\,k\in A$. \ If $m\geq M+s$,
then
\[
x+m\,k=[x+s\,k]+[(m-s)k],
\]
so that $x+m\,k\in A$. \ Thus all sufficiently large $y$ for which \ $y\equiv
r$ ($\operatorname{mod}k$) are in $A$. \ Similarly, for each $c\geq1$, all
sufficiently large $y$ such that \ $y\equiv c\,r$ ($\operatorname{mod}k$) are
in $A$. \ In particular, let $c\geq1$ satisfy%
\[
c\,r\equiv\gcd(r,k)\text{ (}\operatorname{mod}k\text{).}%
\]
Therefore, it follows that all sufficiently large multiples of $\gcd(r,k)$ are
in $A$. \ But $\gcd(r,k)\leq r<k$, which contradicts the minimality of $k$.

\section{Closing Words}

A set $S$ is a \textbf{sum-free set} if none of its sums are in $S$. \ The
structure of such sets is far more complicated than that for sum sets
\cite{C1-mltsum, C2-mltsum, C3-mltsum}. \ A\ simple necessary condition for
$S$ to enjoy regularity is known, but numerical evidence suggests that the
condition fails to be sufficient.

A set $S$ is a \textbf{multisum-free set} if none of its multisums are in $S$.
\ The presence of both unisums \& non-sums in such sets will (almost
certainly) further convolute matters. \ No one has yet studied these, to the
best of our knowledge.

\section{Acknowledgement}

James H. Schmerl was so kind to send me a handwritten letter 30 years ago
outlining his proof \cite{S1-mltsum}. \ I\ am thankful to him for this, as
well as other papers \cite{S2-mltsum, S3-mltsum} relevant to my research at
the time.

\end{document}